# Asymptotic optimality of a cross-validatory predictive approach to linear model selection

**Arijit Chakrabarti[1] and Tapas Samanta[1]**

*Indian Statistical Institute*

**Abstract:** In this article we study the asymptotic predictive optimality of a model selection criterion based on the cross-validatory predictive density, already available in the literature. For a dependent variable and associated explanatory variables, we consider a class of linear models as approximations to the true regression function. One selects a model among these using the criterion under study and predicts a future replicate of the dependent variable by an optimal predictor under the chosen model. We show that for squared error prediction loss, this scheme of prediction performs asymptotically as well as an oracle, where the oracle here refers to a model selection rule which minimizes this loss if the true regression were known.

## Contents

## 1. Introduction

The ultimate goal of modeling in any scientific or sociological investigation is to discover the underlying regular pattern or phenomenon, if any, which controls the data generating mechanism. Although it is almost impossible to imagine that a single model or combinations of a handful will fully capture the intricate functioning of nature or sociological issues, one can always hope to be able to come close. Given a choice of several models and a set of data, a popular method is to choose the model which explains or fits the given data best (in some well-defined sense). However, it is of prime importance that any model that is chosen should be able to predict future observations from the same experiment or process reasonably well and that it does not merely fit the observed data. This is the purpose of predictive model selection.

---

[1] Applied Statistics Division, Indian Statistical Institute, 203 B.T. Road, Kolkata-700108, India, e-mail: arc@isical.ac.in; tapas@isical.ac.in

*AMS 2000 subject classifications:* Primary 62J05; secondary 62F15.

*Keywords and phrases:* cross-validation, oracle, predictive density.

One of the most prominent approaches to predictive model selection is cross-validation (see [17]) and variants thereof. As the name cross-validation suggests, parameters of the population are estimated under each model by using a part of the data (the "estimation set"), while the rest of the data (the "validation set") are predicted using the estimates based on the first group. This is done repeatedly by using "validation sets" comprising different parts of the data, e.g., the whole data could simply be divided into 10 disjoint parts, each part consisting of an equal number of observations and predicted using the rest. If, for a particular model, such predictions match best with the actually observed values, i.e., if the average prediction error is the smallest for it among all the candidate models, it is selected. Optimality properties of classical cross-validatory techniques have been studied, e.g., in [12] and [16].

In the Bayesian literature, several approaches to model selection have been studied with the predictive aspect in mind; see, e.g., [1, 4, 5, 8, 9, 10, 13, 14]. The purpose of this paper is to study the predictive properties of a model selection criterion (see (1.2) below) based on the average of the (log) cross-validatory predictive densities (see (1.1) below) and already available in the literature. Different types of averages (e.g. arithmetic mean, (log) geometric mean) of cross-validatory predictive densities have been studied by several authors ([2], [3], [5], [9] and [14]). Chakrabarti and Ghosh [5] considered an average with respect to disjoint validation sets and studied what should be the optimal proportion of the sample kept for validation in large sample sizes, for the selection of a model closest to the true model (in terms of Kullback-Leibler divergence), and for the selection of the more parsimonious model if two models are equidistant from the truth. Using squared error prediction loss, we show that model selection using criterion (1.2) has an optimality property in predicting a future replicate of the dependent variable (for fixed values of the independent variables), when the true regression is being approximated by a class of candidate linear models. The proofs of the optimality results partly use some general techniques of Li [12] which were later adopted in [16].

In the Bayesian setup, the ordinary predictive density under a model is defined as the integral of the likelihood function of the observed data with respect to the prior distribution of the parameters under the model. Between two competing models, the one having a larger predictive density for the given data seems to be the more appropriate description of the unknown data generating process. In non-subjective Bayesian analysis, it is common to use noninformative priors for the parameters which are typically improper and defined only up to unknown multiplicative constants. In such situations, use of the ordinary predictive density as a model selection criterion will be inappropriate. To get rid of this difficulty, one updates the improper prior by getting a proper posterior based on part of the data (called the training sample) and then integrates the likelihood function of the rest of the data with respect to this posterior, thus giving the cross-validatory predictive density. This is like getting the predictive distribution of part of the data using information obtained from the rest of it. This method of obtaining a cross-validatory predictive density can also be used when one puts a proper prior on the parameters of the model. The cross-validatory predictive density can then be used to get pseudo-Bayes Factors, after appropriate averaging with respect to the different possible choices of the training sample. This line of thought owes its origin to Geisser [7] and Geisser and Eddy [8] and came to prominence through what are referred to as partial Bayes Factors or Intrinsic Bayes Factors ([2], [3], [9], [11] and [15]).

In the next few paragraphs, we describe our setup and the model selection criterion we study. We follow the notations of Shao [16].

Let $\boldsymbol{y}_n = (y_1, \ldots, y_n)'$ be a vector of observations on the dependent (response) variable and let $\boldsymbol{X}_n = (\boldsymbol{x}_1', \ldots, \boldsymbol{x}_n')'$ be an $n \times p_n$ matrix of explanatory variables (which are potentially responsible for the variability in the $y$'s), with $\boldsymbol{x}_i$ associated with $y_i$. Let $\boldsymbol{\mu}_n$ denote $E(\boldsymbol{y}_n | \boldsymbol{X}_n)$, the (unknown) average value of the response variable given the values of the explanatory variables. We further assume that given $\boldsymbol{X}_n$, $\boldsymbol{e}_n = \boldsymbol{y}_n - \boldsymbol{\mu}_n$ has mean vector $\boldsymbol{0}$ and the components of $e_i$ are independent with common variance $\sigma^2$, which could be known or unknown. We are interested in capturing the functional relationship, if any, between $\boldsymbol{\mu}_n$ and $\boldsymbol{X}_n$ which will be most suitable for predictive purposes. We restrict our search within a class of normal linear models. Our model space, denoted $\mathcal{A}_n$, is indexed by $\alpha$, where each $\alpha$ consists of a subset of size $p_n(\alpha)$ $(1 \le p_n(\alpha) \le p_n)$ of $\{1, 2, \ldots, p_n\}$ and the true mean $\boldsymbol{\mu}_n$ is assumed to be linearly related to the corresponding explanatory variables. More specifically, under model $\alpha \in \mathcal{A}_n$, $\boldsymbol{y}_n \sim N(\boldsymbol{\mu}_n(\alpha) = X_n(\alpha)\boldsymbol{\beta}_n(\alpha), \sigma^2 I_n)$ where $\boldsymbol{X}_n(\alpha)$ is the submatrix of $\boldsymbol{X}_n$ consisting of the $p_n(\alpha)$ columns specified by $\alpha$ and $\boldsymbol{\beta}_n(\alpha) \in \Re^{p_n(\alpha)}$. A Bayesian puts a prior on the unknown parameters within each model. We consider standard non-subjective priors (see e.g., [1]) given by

$$\pi_\alpha(\boldsymbol{\beta}_n(\alpha)) \;\; \propto \;\; 1 \qquad \text{if } \sigma^2 \text{ is known, and}$$

$$\pi_\alpha(\boldsymbol{\beta}_n(\alpha), \sigma^2) \;\; \propto \;\; \frac{1}{\sigma^2} \quad \text{if } \sigma^2 \text{ is unknown.}$$

Consider, for example, the case with $\sigma^2$ unknown. Let $\pi_\alpha((\boldsymbol{\beta}_n(\alpha), \sigma^2) | y_{k+1}, \ldots, y_n)$ denote the posterior distribution of the parameters under the model given the observations $(y_{k+1}, \ldots, y_n)$. The cross-validatory predictive density of $(y_1, \ldots, y_k)$ given $(y_{k+1}, \ldots, y_n)$ under model $\alpha$, denoted by the expression $f_\alpha(y_1, \ldots, y_k | y_{k+1}, \ldots, y_n)$, is given by

$$(1.1) \qquad \int f_{\boldsymbol{\beta}_n(\alpha), \sigma^2}(y_1, \ldots, y_k) \pi_\alpha((\boldsymbol{\beta}_n(\alpha), \sigma^2) | y_{k+1}, \ldots, y_n) \, d\boldsymbol{\beta}_n(\alpha) \, d\sigma^2,$$

where $f_{\boldsymbol{\beta}_n(\alpha), \sigma^2}(y_1, \ldots, y_k)$ denotes the density of the $k$ dimensional normal vector, with mean vector given by the first $k$ components of $\boldsymbol{\mu}_n(\alpha)$ and variance-covariance matrix $\sigma^2 I_k$, evaluated at $y_1, \ldots, y_k$. Similarly, the predictive density of any subset $(y_{t_1}, \ldots, y_{t_k})$ of $\boldsymbol{y}$, given the rest of the components of $\boldsymbol{y}$ under this model can be calculated, where $(t_1, \ldots, t_k)$ denotes a subset of $(1, \ldots, n)$. Since a good criterion should not depend too much on the choice of the training sample, we consider the geometric mean of the cross-validatory predictive densities thus obtained by varying the choice of the training sample. The ratio of such geometric means for two models is precisely the Geometric Intrinsic Bayes Factor ([2], [3]). For model $\alpha$, the criterion which we intend to study equals the logarithm of this geometric mean. Thus if we consider a total of $r$ training samples, this logarithm is given by

$$(1.2) \qquad \mathrm{CV}(\alpha) = \frac{1}{r} \sum_{i=1}^{r} \log f_\alpha(y_{t_{1i}}, \ldots, y_{t_{ki}} | \{y_t : t \notin (t_{1i}, \ldots, t_{ki})\}),$$

where $(y_{t_{1i}}, \ldots, y_{t_{ki}})$ is the set of $y$ observations *not* included in the $i$-th training sample. One selects the model $\hat\alpha_n \in \mathcal{A}_n$ which maximizes $\mathrm{CV}(\alpha)$.

Once a model is thus selected, we use the mean of the predictive distribution of $\boldsymbol{y}_n^{\mathrm{new}}$, given the observed $\boldsymbol{y}_n$ under the selected model, as the predictor for a future replicate $\boldsymbol{y}_n^{\mathrm{new}}$ of the response variable for the same value $\boldsymbol{X}_n$ of the explanatory variables. An easy calculation shows that this turns out to be the least squares estimate $X(\hat\alpha_n)\hat{\boldsymbol{\beta}}_n(\hat\alpha_n)$ where $\hat{\boldsymbol{\beta}}_n(\alpha) = P_n(\alpha)\boldsymbol{y}_n$ and $P_n(\alpha) = X_n(\alpha)(X_n(\alpha)'X_n(\alpha))^{-1}X_n(\alpha)$ is the usual projection matrix.

Our goal is to evaluate this prediction scheme under the true regression using squared error prediction loss. Under the true $\boldsymbol{\mu}_n$, the future replicate $\boldsymbol{y}_n^{\text{new}}$ will be independent of the original observations $\boldsymbol{y}_n$. The quality of any predictor $\delta(\boldsymbol{y}_n)$ of $\boldsymbol{y}_n^{\text{new}}$ based on $\boldsymbol{y}_n$ can be evaluated by the average prediction error $E_{\boldsymbol{\mu}_n}(\frac{1}{n}||\boldsymbol{y}_n^{\text{new}} - \delta(\boldsymbol{y}_n)||^2)$, where $E_{\boldsymbol{\mu}_n}$ denotes expectation with respect to the joint distribution of $(\boldsymbol{y}_n^{\text{new}}, \boldsymbol{y}_n)$ when $\boldsymbol{\mu}_n$ is the true unknown mean. This expectation will be small if, for any fixed $\boldsymbol{y}_n$, $E_{\boldsymbol{\mu}_n}(\frac{1}{n}||\boldsymbol{y}_n^{\text{new}} - \delta(\boldsymbol{y}_n)||^2|\boldsymbol{y}_n)$ is also small. As observed before, the predictor $\delta(\boldsymbol{y}_n)$ we want to evaluate is the same as the least squares predictive estimate of $\boldsymbol{y}_n^{\text{new}}$ under the chosen model $\hat{\alpha}_n$. Now note that for any given fixed model $\alpha$, the least squares predictive estimate is given by $\delta(\boldsymbol{y}_n) = \delta(\boldsymbol{y}_n)(\alpha) = \hat{\boldsymbol{\mu}}_n(\alpha) = X_n(\alpha)\hat{\boldsymbol{\beta}}_n(\alpha)$. A simple algebra shows that the above conditional expectation is, up to a constant which does not depend on $\alpha$, equal to

$$(1.3) \qquad L_n(\alpha) = \frac{||\boldsymbol{\mu}_n - \hat{\boldsymbol{\mu}}_n(\alpha)||^2}{n}.$$

Hence the conditional expectation will be minimized for a certain $\alpha$ if $L_n(\alpha)$ is minimized. If we knew the true $\boldsymbol{\mu}_n$, we could find the model which minimizes this $L_n(\alpha)$ for each $\boldsymbol{y}_n$. We shall call this the oracle model, denoted $\alpha_n^L$. The best any procedure can achieve is to do as well as the oracle in the limit in terms of the loss as the sample size grows to infinity.

We show in the following sections of this article that under certain conditions, minimizing $\text{CV}(\alpha)$ with respect to $\alpha$ is asymptotically equivalent to minimizing $L_n(\alpha)$. Using this fact it is shown that the ratio of $L_n(\alpha_n^L)$ to $L_n(\hat{\alpha}_n)$ tends to 1 in probability, whereby establishing the optimum asymptotic behavior of criterion (1.2) in the problem of prediction of a set of future observations.

In Sections 2 and 3 we consider the case where the true model is not in the model space – the proposed models are only approximations to the truth. In Section 2 we consider the case when $\sigma^2$ is known. We show that under certain assumptions, the model selection procedure under study performs as well as the oracle asymptotically in the sense that the ratio of their losses tends to one in probability. In Section 3, we consider the more realistic situation when $\sigma^2$ is unknown. Under appropriate conditions it is shown that this procedure also achieves the oracle asymptotically in this case. As a validation of this method, we next consider in Section 4 the question of whether, under the assumption that the true model is indeed included in the model space, we do equally well in terms of hitting the oracle loss asymptotically. It is shown that this model selection procedure chooses the correct model with smallest dimension with probability tending to one in addition to being asymptotically optimal in terms of hitting the oracle. Some concluding remarks are made in Section 5. Technical proofs of most of the results are given in the Appendix.

For notational simplicity we write $\boldsymbol{y}, \boldsymbol{\mu}, \boldsymbol{e}, X(\alpha), \boldsymbol{\beta}(\alpha)$ and $P(\alpha)$ in place of $\boldsymbol{y}_n$, $\boldsymbol{\mu}_n, \boldsymbol{e}_n, X_n(\alpha), \boldsymbol{\beta}_n(\alpha)$ and $P_n(\alpha)$ respectively, dropping the suffix $n$ for the rest of the paper.

## 2. Basic results – case with $\boldsymbol{\sigma^2}$ known

In this section we take the "model false" point of view that the models are only approximations to the truth but none of them is actually true. We show that under certain conditions, the model selection procedure under study is asymptotically optimal in the sense of performing as well as the oracle defined above.

As described in the introduction, the model selection criterion under consideration is an average of the cross-validatory predictive density

$$f_\alpha(y_1, \ldots, y_k | y_{k+1}, \ldots, y_n)$$

under model $\alpha$, over suitable choices of the "training sample" $\{y_{k+1}, \ldots, y_n\}$. *We do not recommend here any particular choice of the training samples; our results hold as long as each $y_i$, $1 \le i \le n$, appears in the same number of training samples chosen (which will be assumed throughout the paper).*

Let $\boldsymbol{y}_i$, $i = 1, \ldots, r$ be the $r$ training samples (each of size $n-k$). For each $\boldsymbol{y}_i$, let $\boldsymbol{\mu}_i$ and $\boldsymbol{e}_i$ be the subvector of $\boldsymbol{\mu}$ and $\boldsymbol{e}$ corresponding to the labels of the components of $\boldsymbol{y}_i$ and $X_i(\alpha)$ be the submatrix of $X(\alpha)$ consisting of the corresponding rows of it. Also, let $\hat{\boldsymbol{\beta}}_i(\alpha) = [X_i'(\alpha)X_i(\alpha)]^{-1}X_i'(\alpha)\boldsymbol{y}_i$, $P_i(\alpha) = X_i(\alpha)[X_i'(\alpha)X_i(\alpha)]^{-1}X_i'(\alpha)$, $i = 1, \ldots, r$. It will be assumed throughout that $(n-k) \to \infty$ and $X_i'(\alpha)X_i(\alpha)$ is nonsingular for each $i$ and $\alpha$. With the standard non-subjective prior $\pi(\boldsymbol{\beta}(\alpha)) = $ constant, we have a closed form expression for the cross-validatory predictive density. An alternative equivalent criterion, which is to be minimized with respect to $\alpha$, is

$$
\begin{aligned}
\Gamma(\alpha) \;=\; & \frac{1}{n}(\boldsymbol{y} - X(\alpha)\hat{\boldsymbol{\beta}}(\alpha))'(\boldsymbol{y} - X(\alpha)\hat{\boldsymbol{\beta}}(\alpha)) \\
& - \frac{1}{r}\sum_{i=1}^{r}\frac{1}{n}(\boldsymbol{y}_i - X_i(\alpha)\hat{\boldsymbol{\beta}}_i(\alpha))'(\boldsymbol{y}_i - X_i(\alpha)\hat{\boldsymbol{\beta}}_i(\alpha)) \\
& + \frac{1}{r}\sum_{i=1}^{r}\frac{\sigma^2}{n}\log\left(\frac{|X'(\alpha)X(\alpha)|}{|X_i'(\alpha)X_i(\alpha)|}\right).
\end{aligned}
$$
(2.1)

Note that $\Gamma(\alpha)$ is equal to the negative of the criterion (1.2) up to an additive constant. We will prove that minimization of $\Gamma(\alpha)$ is equivalent to minimization of the loss $L_n(\alpha)$ (defined in (1.3)) in an appropriate asymptotic sense and this will lead to the desired asymptotic (predictive) optimality of the criterion under consideration.

Note that the loss $L_n(\alpha)$ defined in (1.3) can be written as

$$nL_n(\alpha) = n\Delta_n(\alpha) + \boldsymbol{e}'P(\alpha)\boldsymbol{e}$$

where $n\Delta_n(\alpha) = \boldsymbol{\mu}'(I - P(\alpha))\boldsymbol{\mu}$ and let

$$nR_n(\alpha) = E(nL_n(\alpha)) = n\Delta_n(\alpha) + \sigma^2 p_n(\alpha).$$

One of the key assumptions under which we prove our results is the following condition ([12], [16]):

$$\sum_{\alpha \in \mathcal{A}_n} \frac{1}{[nR_n(\alpha)]^m} \to 0$$
(2.2)

for some positive integer $m$ for which $E(e_1^{4m}) < \infty$. We also assume

$$\frac{p_n \lambda_n}{\min_{\alpha \in \mathcal{A}_n} nR_n(\alpha)} \to 0,$$
(2.3)

where $\lambda_n = \log(n/(n-k))$.

For certain remarks justifying these assumptions, see [12] and [16]. In particular, it is argued in these papers using several concrete examples, that condition (2.2) is a natural one when the dimension $p_n$ of the largest model grows with sample size. Also, if $p_n$ remains bounded, $nR_n(\alpha)$ is expected to go to $\infty$ for all $\alpha$ as the sample size increases, if the candidate models are separated from the truth. That $\min_\alpha nR_n(\alpha) \to \infty$ is assumption A.3$'$ of Li [12] and as remarked therein, it is a quite reasonable assumption if $p_n$ grows with $n$. Condition (2.3) requires that $\min_\alpha nR_n(\alpha) \to \infty$ at a suitable rate. Under condition (3.3) below ([16], condition (2.5)), (2.3) holds if $(p_n\lambda_n)/n \to 0$.

It is important to note that we also need to assume $(n-k)/n \to 0$ to prove our results (see e.g. (6.10)). This addresses an important question about the required size of the training sample. We, however, do not claim that it is a necessary condition for asymptotic predictive optimality.

We now consider the criterion $\Gamma(\alpha)$ as defined in (2.1). Since $X(\alpha)\hat{\boldsymbol{\beta}}(\alpha) = P(\alpha)\boldsymbol{y}$,

$$
\begin{aligned}
&\frac{1}{n}(\boldsymbol{y} - X(\alpha)\hat{\boldsymbol{\beta}}(\alpha))'(\boldsymbol{y} - X(\alpha)\hat{\boldsymbol{\beta}}(\alpha)) \\
=\ &\frac{1}{n}\boldsymbol{y}'(I - P(\alpha))\boldsymbol{y} \\
=\ &\frac{1}{n}\boldsymbol{e}'\boldsymbol{e} + L_n(\alpha) - \frac{2}{n}\boldsymbol{e}'P(\alpha)\boldsymbol{e} + \frac{2}{n}\boldsymbol{e}'(I - P(\alpha))\boldsymbol{\mu}.
\end{aligned}
$$

(2.4)

Similarly,

$$
\begin{aligned}
&\frac{1}{r}\sum_{i=1}^r \frac{1}{n}(\boldsymbol{y}_i - X_i(\alpha)\hat{\boldsymbol{\beta}}_i(\alpha))'(\boldsymbol{y}_i - X_i(\alpha)\hat{\boldsymbol{\beta}}_i(\alpha)) \\
=\ &\frac{n-k}{n^2}\boldsymbol{e}'\boldsymbol{e} + \frac{1}{nr}\sum_{i=1}^r \boldsymbol{\mu}_i'(I - P_i(\alpha))\boldsymbol{\mu}_i - \frac{1}{nr}\sum_{i=1}^r \boldsymbol{e}_i'P_i(\alpha)\boldsymbol{e}_i \\
&+ \frac{2}{nr}\sum_{i=1}^r \boldsymbol{e}_i'(I - P_i(\alpha))\boldsymbol{\mu}_i.
\end{aligned}
$$

(2.5)

We first state two auxiliary results.

**Lemma 2.1.** *Under conditions (2.2) and (2.3),*

$$
\frac{1}{n}(\boldsymbol{y} - X(\alpha)\hat{\boldsymbol{\beta}}(\alpha))'(\boldsymbol{y} - X(\alpha)\hat{\boldsymbol{\beta}}(\alpha)) = \frac{1}{n}\boldsymbol{e}'\boldsymbol{e} + L_n(\alpha) + o_p(L_n(\alpha))
$$

*uniformly in $\alpha \in \mathcal{A}_n$.*

By saying $Z_n(\alpha) = o_p(L_n(\alpha))$ uniformly in $\alpha$, we mean $\max_\alpha |Z_n(\alpha)|/L_n(\alpha) \xrightarrow{p} 0$.

**Lemma 2.2.** *Suppose that conditions (2.2) and (2.3) hold and $(n-k)/n \to 0$. Then*

$$
\frac{1}{r}\sum_{i=1}^r \frac{1}{n}(\boldsymbol{y}_i - X_i(\alpha)\hat{\boldsymbol{\beta}}_i(\alpha))'(\boldsymbol{y}_i - X_i(\alpha)\hat{\boldsymbol{\beta}}_i(\alpha)) = \frac{n-k}{n^2}\boldsymbol{e}'\boldsymbol{e} + o_p(L_n(\alpha)),
$$

*uniformly in $\alpha \in \mathcal{A}_n$.*

Proofs of Lemma 2.1 and Lemma 2.2 are given in the Appendix.

In order to prove the main result of this section we need to assume another condition which is given below.

Let

$$(2.6) \qquad a_{in}(\alpha) = \log \left\{ \frac{(n-k)^{p_n(\alpha)} |X'(\alpha) X(\alpha)|}{n^{p_n(\alpha)} |X'_i(\alpha) X_i(\alpha)|} \right\}.$$

We assume

$$(2.7) \qquad \max_{\alpha \in \mathcal{A}_n} \frac{\frac{1}{r} \sum_{i=1}^{r} a_{in}(\alpha)}{n R_n(\alpha)} \to 0.$$

**Remark 2.1.** Let $\boldsymbol{x}'_1(\alpha), \ldots, \boldsymbol{x}'_n(\alpha)$ be the $n$ rows of $X(\alpha)$. If these $n$ rows are "similar", e.g., if they can be thought of as (independent) realizations of a random vector $\boldsymbol{x}$ and $p_n$ is small compared to both $n-k$ and $n$, then

$$\left| \frac{X'(\alpha) X(\alpha)}{n} \right| = \left| \frac{1}{n} \sum_{j=1}^{n} \boldsymbol{x}_j(\alpha) \boldsymbol{x}'_j(\alpha) \right| \approx |E(\boldsymbol{x}\boldsymbol{x}')|$$

$$\text{and similarly } \left| \frac{X'_i(\alpha) X_i(\alpha)}{n-k} \right| \approx |E(\boldsymbol{x}\boldsymbol{x}')|.$$

In this case, it follows that $a_{in}(\alpha) \approx 0$. In such a situation, assumption (2.7) seems to be quite reasonable.

Now note that (2.3) and (2.7) will imply that the third term in the right hand side of (2.1) is also of the order $o_p(L_n(\alpha))$ uniformly in $\alpha \in \mathcal{A}_n$. Thus

$$\Gamma(\alpha) = \text{constant} + L_n(\alpha) + o_p(L_n(\alpha)) \text{ uniformly in } \alpha \in \mathcal{A}_n$$

which implies minimization of $\Gamma(\alpha)$ is essentially equivalent to minimization of $L_n(\alpha)$ in an appropriate asymptotic sense and we have the following result.

**Theorem 2.1.** *Suppose that conditions (2.2), (2.3) and (2.7) hold and $(n-k)/n \to 0$. Then we have the following results.*

(a) $\Gamma(\alpha) = \frac{k}{n^2} \boldsymbol{e}' \boldsymbol{e} + L_n(\alpha) + o_p(L_n(\alpha))$ *uniformly in $\alpha \in \mathcal{A}_n$.*

(b) *The model selection rule under study is asymptotically optimal in the sense that*

$$\frac{L_n(\hat{\alpha}_n)}{\min_{\alpha \in \mathcal{A}_n} L_n(\alpha)} \xrightarrow{p} 1$$

*where $\hat{\alpha}_n$ is as defined in Section 1.*

Proof of Theorem 2.1 is given in the Appendix.

## 3. Case with $\sigma^2$ unknown

We now consider the more realistic situation when the variance $\sigma^2$ is unknown. The standard non-subjective prior in this case is $\pi(\boldsymbol{\beta}(\alpha), \sigma^2) \propto \frac{1}{\sigma^2}$ under model $\alpha$. Interestingly, the results in this case follow from the basic results obtained in Section 2. We consider here the ("model false") setup and assumptions of Section 2.

Let $\boldsymbol{y}_i, i = 1, \ldots, r$ be the $r$ training samples chosen. The cross-validatory predictive density under model $\alpha$ for a training sample $\boldsymbol{y}_i$ is given by

$$\frac{|X_i'(\alpha)X_i(\alpha)|^{\frac{1}{2}}}{|X'(\alpha)X(\alpha)|^{\frac{1}{2}}} \times \frac{[(\boldsymbol{y} - X(\alpha)\hat{\boldsymbol{\beta}}(\alpha))'(\boldsymbol{y} - X(\alpha)\hat{\boldsymbol{\beta}}(\alpha))]^{-\frac{n}{2}}}{[(\boldsymbol{y}_i - X_i(\alpha)\hat{\boldsymbol{\beta}}_i(\alpha))'(\boldsymbol{y}_i - X_i(\alpha)\hat{\boldsymbol{\beta}}_i(\alpha))]^{-\frac{n-k}{2}}}$$

up to a multiplicative constant.

Our criterion (to be minimized with respect to $\alpha$), which is an average over the $r$ training samples, is given by

$$(3.1) \qquad \Gamma(\alpha) = \log[S(\alpha)] - \frac{n-k}{nr}\sum_{i=1}^{r}\log[S_i(\alpha)] + \frac{1}{nr}\sum_{i=1}^{r}\log\frac{|X'(\alpha)X(\alpha)|}{|X_i'(\alpha)X_i(\alpha)|}$$

where $S(\alpha) = (\boldsymbol{y} - X(\alpha)\hat{\boldsymbol{\beta}}(\alpha))'(\boldsymbol{y} - X(\alpha)\hat{\boldsymbol{\beta}}(\alpha))$ and $S_i(\alpha) = (\boldsymbol{y}_i - X_i(\alpha)\hat{\boldsymbol{\beta}}_i(\alpha))'(\boldsymbol{y}_i - X_i(\alpha)\hat{\boldsymbol{\beta}}_i(\alpha))$.

Note that $\Gamma(\alpha) = (k/n)\log(n\sigma^2) + \Gamma_1(\alpha)$ where

$$(3.2) \quad \Gamma_1(\alpha) = \log\left[\frac{S(\alpha)}{n\sigma^2}\right] - \frac{n-k}{nr}\sum_{i=1}^{r}\log\left[\frac{S_i(\alpha)}{n\sigma^2}\right] + \frac{1}{nr}\sum_{i=1}^{r}a_{in}(\alpha) + \frac{1}{n}p_n(\alpha)\lambda_n,$$

$a_{in}(\alpha)$ is as defined in (2.6) and $\lambda_n = \log(n/(n-k))$. Therefore, minimizing $\Gamma(\alpha)$ (with respect to $\alpha$) is equivalent to minimizing $\Gamma_1(\alpha)$ for all $\sigma$. Let

$$u_n(\alpha) = \log\left[\frac{\boldsymbol{e}'\boldsymbol{e}}{n\sigma^2} + \frac{1}{\sigma^2}L_n(\alpha)\right].$$

In order to prove the asymptotic optimality of this method, we first note in Lemma 3.1 below that $\Gamma_1(\alpha)$ is asymptotically equivalent to $u_n(\alpha)$ and this in turn implies the desired conclusion as stated in Theorem 3.1. We prove these results by invoking certain conditions which we describe below.

We first make the following assumption (see [16], condition (2.5)):

$$(3.3) \qquad \liminf_{n \to \infty}\min_{\alpha}\Delta_n(\alpha) > 0$$

where $\Delta_n(\alpha)$ is as defined in Section 2. This may be thought of as an identifiability condition on the models in the model space, as appears in the discussion of Mervyn Stone on [16]. We further assume that

$$(3.4) \qquad \frac{n-k}{n}\log n \to 0, \quad \frac{p_n\lambda_n}{n} \to 0 \text{ and } \frac{1}{n}\sum_{i=1}^{n}\mu_i^2 \text{ is bounded},$$

$$(3.5) \qquad \frac{1}{nr}\sum_{i=1}^{r}a_{in}(\alpha) \to 0,$$

and

$$(3.6) \qquad \sum_{i=1}^{r}\log(S_i) > 0$$

with probability tending to 1, where $S_i$ is equal to $S_i(\alpha)$ with $\alpha$ as the full model, i.e., $\alpha = \{1, \ldots, p_n\}$. One can give sufficient conditions for (3.6) based on the relative magnitude of $r$ and $(n-k)$ as $n \to \infty$, to the effect that $r$ is not too large compared with $n-k$ which is the case for most practically implementable schemes. We, however, do not record the details here. The final results of this section are now stated below.

**Lemma 3.1.** *Under conditions (3.3)–(3.6),*

$$(3.7) \qquad \Gamma_1(\alpha) = u_n(\alpha) + o_p(u_n(\alpha)) \ \text{uniformly in } \alpha.$$

**Theorem 3.1.** *Under conditions (3.3)–(3.6),*

$$(3.8) \qquad \frac{L_n(\hat{\alpha}_n)}{L_n(\alpha_n^L)} \xrightarrow{p} 1.$$

Both Lemma 3.1 and Theorem 3.1 are proved in the Appendix.

## 4. The "model true" case and consistency

We now show that if some model in the model space is true, the model selection procedure under study chooses the correct model of the smallest dimension in addition to being asymptotically optimal. Thus this procedure not only captures the truth but at the same time is as parsimonious as possible. Although the assumption of a true model may not seem to be very realistic, our result in this section provides a validation of the method. We, however, consider only the simpler case when $\sigma^2$ is known.

As in [16], let $\mathcal{A}_n^c \subset \mathcal{A}_n$ denote all the proposed models that are actually correct. Thus for $\alpha \in \mathcal{A}_n^c$, $\boldsymbol{\mu} = X(\alpha)\boldsymbol{\beta}(\alpha)$ for some $\boldsymbol{\beta}(\alpha) \in \Re^{p_n(\alpha)}$. In Section 2 we assumed that $\mathcal{A}_n^c$ is empty. It is important to note that all the results of Section 2 with $\mathcal{A}_n$ replaced by $\mathcal{A}_n - \mathcal{A}_n^c$ hold under the corresponding assumptions with $\mathcal{A}_n$ replaced by $\mathcal{A}_n - \mathcal{A}_n^c$. In particular, if

$$(4.1) \qquad \sum_{\alpha \in \mathcal{A}_n - \mathcal{A}_n^c} \frac{1}{[nR_n(\alpha)]^m} \to 0$$

for some positive integer $m$ for which $E(e_1^{4m}) < \infty$ and

$$(4.2) \qquad \frac{p_n \lambda_n}{\min_{\alpha \in \mathcal{A}_n - \mathcal{A}_n^c} nR_n(\alpha)} \to 0,$$

with $\lambda_n = \log(n/(n-k))$, then

$$(4.3) \qquad \Gamma(\alpha) = \frac{k}{n^2} \boldsymbol{e}' \boldsymbol{e} + L_n(\alpha) + o_p(L_n(\alpha))$$

uniformly in $\alpha \in \mathcal{A}_n - \mathcal{A}_n^c$.

For $\alpha \in \mathcal{A}_n^c$, $(I - P(\alpha))\boldsymbol{\mu} = \boldsymbol{0}$ and $(I - P_i(\alpha))\boldsymbol{\mu}_i = 0 \ \forall \ i$. Therefore, from (2.1), (2.4) and (2.5) we have for $\alpha \in \mathcal{A}_n^c$

$$(4.4) \quad \Gamma(\alpha) = \frac{k}{n^2} \boldsymbol{e}' \boldsymbol{e} - \frac{1}{n} \boldsymbol{e}' P(\alpha) \boldsymbol{e} + \frac{1}{nr} \sum_{i=1}^{r} \boldsymbol{e}_i' P_i(\alpha) \boldsymbol{e}_i + \frac{\sigma^2}{nr} \sum_{i=1}^{r} \log\left(\frac{|X'(\alpha)X(\alpha)|}{|X_i'(\alpha)X_i(\alpha)|}\right).$$

Also $L_n(\alpha) = \frac{1}{n} \boldsymbol{e}' P(\alpha) \boldsymbol{e}$ for $\alpha \in \mathcal{A}_n^c$.

We now assume that

$$(4.5) \qquad \limsup_{n \to \infty} \sum_{\alpha \in \mathcal{A}_n^c} \frac{1}{[p_n(\alpha)]^m} < \infty.$$

for some positive integer $m$ such that $E(e_1^{4m}) < \infty$ (condition (3.10) of Shao [16]), and

$$(4.6) \qquad \max_{\alpha \in \mathcal{A}_n^c} \frac{\frac{1}{r} \sum_{i=1}^{r} a_{in}(\alpha)}{p_n(\alpha)\lambda_n} \to 0$$

with $\lambda_n = \log(\frac{n}{n-k})$ and $a_{in}(\alpha)$ as defined in (2.6). See Remark 2.1 in this context. Let $\alpha_n^c$ be the model $\alpha$ in $\mathcal{A}_n^c$ with smallest dimension. Using the above, we now have

**Proposition 4.1.** *Under conditions (4.1), (4.2), (4.5) and (4.6)*

$$(4.7) \qquad \Gamma(\alpha) = \frac{k}{n^2}\boldsymbol{e}'\boldsymbol{e} + \frac{1}{n}\lambda_n\sigma^2 p_n(\alpha) + o_p(\frac{1}{n}\lambda_n\sigma^2 p_n(\alpha))$$

*uniformly in $\alpha \in \mathcal{A}_n^c$, and*

$$(4.8) \qquad \Gamma(\alpha_n^c) = \frac{k}{n^2}\boldsymbol{e}'\boldsymbol{e} + o_p(L_n(\alpha)) \text{ uniformly in } \alpha \in \mathcal{A}_n - \mathcal{A}_n^c.$$

Proof of Proposition 4.1 is given in the Appendix.

Keeping in mind the above facts, we now proceed towards proving that this model selection rule chooses the most parsimonious correct model as claimed in Theorem 4.1 below. Towards this we first observe that (4.3) and (4.8) imply

$$\max_{\alpha \in \mathcal{A}_n - \mathcal{A}_n^c} (\Gamma(\alpha_n^c) - \frac{k}{n^2}\boldsymbol{e}'\boldsymbol{e})/(\Gamma(\alpha) - \frac{k}{n^2}\boldsymbol{e}'\boldsymbol{e}) < 1$$

with probability tending to 1. It then follows that

$$(4.9) \qquad P[\Gamma(\alpha_n^c) \leq \Gamma(\alpha) \; \forall \alpha \in \mathcal{A}_n - \mathcal{A}_n^c] \to 1.$$

We now try to find some conditions under which

$$(4.10) \qquad P[\Gamma(\alpha_n^c) \leq \Gamma(\alpha) \; \forall \alpha \in \mathcal{A}_n^c] \to 1.$$

Let $n[\Gamma(\alpha) - \Gamma(\alpha_n^c)] = Z_n(\alpha)$. It is enough to show that

$$(4.11) \qquad P[Z_n(\alpha) \geq 0 \; \forall \alpha \in \mathcal{A}_n^c] \to 1.$$

Now,

$$P[Z_n(\alpha) < 0 \text{ for some } \alpha \in \mathcal{A}_n^c]$$

$$\leq \sum_{\alpha \in \mathcal{A}_n^c} P[Z_n(\alpha) < 0]$$

$$\leq \sum_{\alpha \in \mathcal{A}_n^c} P[|Z_n(\alpha) - E(Z_n(\alpha))| > E(Z_n(\alpha))]$$

$$(4.12) \qquad \leq \sum_{\alpha \in \mathcal{A}_n^c} \frac{E|Z_n(\alpha) - E(Z_n(\alpha))|^{2m}}{[E(Z_n(\alpha))]^{2m}}.$$

From (4.4)

$$(4.13) \quad Z_n(\alpha) - E(Z_n(\alpha)) = \frac{1}{r}\sum_{i=1}^{r} \boldsymbol{e}_i'[P_i(\alpha) - P_i(\alpha_n^c)]\boldsymbol{e}_i - \boldsymbol{e}'[P(\alpha) - P(\alpha_n^c)]\boldsymbol{e}$$

and $E(Z_n(\alpha))$ can be written as

$$\frac{1}{\sigma^2} E(Z_n(\alpha)) = [p_n(\alpha) - p_n(\alpha_n^c)]\lambda_n + \frac{1}{r}\sum_{i=1}^{r}[a_{in}(\alpha) - a_{in}(\alpha_n^c)]$$

where $a_{in}(\alpha)$ is as defined in (2.6). If we assume

$$(4.14) \qquad \frac{1}{r}\sum_{i=1}^{r}[a_{in}(\alpha) - a_{in}(\alpha_n^c)] = o_p([p_n(\alpha) - p_n(\alpha_n^c)]\lambda_n)$$

uniformly in $\alpha \in \mathcal{A}_n^c$, then

$$(4.15) \qquad \frac{1}{\sigma^2} E(Z_n(\alpha)) = [p_n(\alpha) - p_n(\alpha_n^c)]\lambda_n + o_p([p_n(\alpha) - p_n(\alpha_n^c)]\lambda_n)$$

uniformly in $\alpha \in \mathcal{A}_n^c$. Noting that $P(\alpha) - P(\alpha_n^c)$ and $P_i(\alpha) - P_i(\alpha_n^c)$ are projection matrices and the first term on the right hand side of (4.13) can be expressed as $\boldsymbol{e}'M\boldsymbol{e}$ for some matrix $M$, and using Theorem 2 of Whittle [18] or inequality (6.2) of the Appendix we have

$$E|Z_n(\alpha) - E(Z_n(\alpha))|^{2m} \leq \mathrm{constant}[p_n(\alpha) - p_n(\alpha_n^c)]^m.$$

It then follows from (4.12) and (4.15) that (4.11) holds if

$$(4.16) \qquad \sum_{\alpha \in \mathcal{A}_n^c} \frac{1}{\lambda_n^{2m}[p_n(\alpha) - p_n(\alpha_n^c)]^m} \to 0.$$

Thus we finally have the following.

**Theorem 4.1.** *Under conditions (4.1), (4.2), (4.5), (4.6), (4.14) and (4.16),*

$$(4.17) \qquad P[\hat{\alpha}_n = \alpha_n^c] \to 1.$$

It is proved in the Appendix that under (4.1) and (4.2)

$$(4.18) \qquad \max_{\alpha \in \mathcal{A}_n - \mathcal{A}_n^c} \frac{L_n(\alpha_n^c)}{L_n(\alpha)} \xrightarrow{p} 0.$$

Since $L_n(\alpha_n^c) \leq L_n(\alpha) \ \forall \alpha \in \mathcal{A}_n^c$, Theorem 4.1 and (4.18) imply the following.

**Theorem 4.2.** *Under the conditions of Theorem 4.1, one has*

$$(4.19) \qquad L_n(\hat{\alpha}_n)/L_n(\alpha_n^L) \xrightarrow{p} 1.$$

## 5. Concluding remarks

In this article we have studied predictive optimality of a cross-validatory Bayesian approach to model selection in the context of selecting from among a set of linear models. It has been shown that this method predicts as well as the oracle as the sample size grows. In addition, it has been shown that in case the space of candidate models contains at least one correct model, this method chooses the correct model with the smallest dimension with probability tending to one as sample size grows. Thus the method has two important facets – one of an optimal predictor and the other of a selection criterion which does not unnecessarily choose a complex model when simpler ones are apt.

Needless to say, this article has not addressed some interesting related issues. First, it will be interesting to see how this method works when it is applied in the setup of generalized linear models, through theoretical investigation and simulation. Another focus of recent research is the case when the number of potential parameters in the models is very large, e.g., when it is of the same order as the number of observations. Asymptotic optimality studies in such setup, even for the normal linear models will be a really challenging task. Also, we have not touched upon the computational aspect of this method, which becomes important if the number of potential regressors and number of models in the model space get large. We, however, emphasize that one rarely considers the set of all $2^p$ possible models if $p$ regressors are available. For example, one can use expert knowledge about the problem under study and start with a pruned list of models or one can take a nested sequence of models (thereby restricting the total number of models to at most $p$). Li ([12], Example 1) considered a situation where the $p$ regressors are arranged in decreasing order of importance. He then considered $p$ models, the $\alpha$-th model consisting of the first $\alpha$ regressors in this ordered arrangement. See in this context Examples 1 and 2 of [16] where the number of models under consideration is fixed although the number of parameters may grow with sample size. Last but not the least, as we commented before, the requirement that $k/n \to 1$ is only a sufficient condition; a careful study of the necessity of this condition is in order. In some examples, we have observed that $k/n \to c$ for any $c \in (0, 1)$ is also sufficient to achieve good optimality results similar to ones we have obtained in this paper. Some theoretical investigations and simulation studies will hopefully prove conclusive to find the optimal $k$. It is worth mentioning that in a related problem Chakrabarti and Ghosh [5] made interesting observations regarding this issue which can be a starting point for such investigation.

## Appendix

We present in this section proofs of some of the results of the earlier sections. We will need bounds for the moments of linear and quadratic forms in $\boldsymbol{e}$. Let $A = (a_{ij})$ be a non-random $n \times n$ matrix and $\boldsymbol{b}$ be a non-random $n$-vector. Then by Theorem 2 of Whittle [18],

$$(6.1) \qquad E(|\boldsymbol{e}'\boldsymbol{b}|^{2m}) \leq C_1(||\boldsymbol{b}||^2)^m, \text{ and}$$

$$(6.2) \qquad E|\boldsymbol{e}'A\boldsymbol{e} - E(\boldsymbol{e}'A\boldsymbol{e})|^{2m} \leq C_2(\sum_i \sum_j a_{ij}^2)^m$$

for some constants $C_1, C_2 > 0$ and for positive integer $m$ for which $E(e_1^{4m}) < \infty$.

Below $\max_\alpha$ will mean maximum over $\alpha \in \mathcal{A}_n$.

*Proof of Lemma 2.1.* As shown in Li ([12], p. 970), using Theorem 2 of Whittle [18] or inequalities (6.1) and (6.2) stated above, and condition (2.2), we have

$$(6.3) \qquad \max_\alpha \frac{|\boldsymbol{e}'P(\alpha)\boldsymbol{e} - \sigma^2 p_n(\alpha)|}{nR_n(\alpha)} \xrightarrow{p} 0, \text{ and}$$

$$(6.4) \qquad \max_\alpha \frac{|\boldsymbol{e}'(I - P(\alpha))\boldsymbol{\mu}|}{nR_n(\alpha)} \xrightarrow{p} 0.$$

Also, from (6.3)

$$(6.5) \qquad \max_\alpha |\frac{L_n(\alpha)}{R_n(\alpha)} - 1| = \max_\alpha \frac{|\boldsymbol{e}'P(\alpha)\boldsymbol{e} - \sigma^2 p_n(\alpha)|}{nR_n(\alpha)} \xrightarrow{p} 0.$$

Lemma 2.1 now follows from (2.3), (2.4), (6.3), (6.4) and (6.5).     □

*Proof of Lemma 2.2.* Let

$$T_1 = \frac{1}{r}\sum_{i=1}^r \boldsymbol{\mu}_i'(I - P_i(\alpha))\boldsymbol{\mu}_i, \ T_2 = \frac{1}{r}\sum_{i=1}^r \boldsymbol{e}_i'P_i(\alpha)\boldsymbol{e}_i \text{ and } T_3 = \frac{1}{r}\sum_{i=1}^r \boldsymbol{e}_i'(I - P_i(\alpha))\boldsymbol{\mu}_i.$$

Then, in view of (2.5), the left hand side of the equality claimed in Lemma 2.2 can be written as

$$\frac{n-k}{n^2}\boldsymbol{e}'\boldsymbol{e} + \frac{1}{n}(T_1 - T_2 + 2T_3).$$

We shall prove that

$$(6.6) \qquad T_j/n = o_p(L_n(\alpha)) \text{ uniformly in } \alpha$$

for $j = 1, 2, 3$.

We fix a training sample $\boldsymbol{y}_1 = (y_1, y_2, \ldots, y_{n-k})'$. Let

$$X(\alpha) = \begin{pmatrix} X_1 \\ X_{1c} \end{pmatrix} \text{ and } I - P(\alpha) = \begin{pmatrix} A \\ B \end{pmatrix}$$

where $X_1$ and $X_{1c}$ are the submatrices consisting of the first $n - k$ rows and the last $k$ rows of $X$, respectively, and $A$ and $B$ are analogous submatrices of $I - P(\alpha)$. Then

$$(6.7) \qquad\qquad \boldsymbol{\mu}'(I - P(\alpha))\boldsymbol{\mu} = \boldsymbol{\mu}'B'B\boldsymbol{\mu} + \boldsymbol{\mu}'A'A\boldsymbol{\mu}, \text{ and}$$
$$(6.8) \qquad \boldsymbol{\mu}'(I - P(\alpha))\boldsymbol{\mu} - \boldsymbol{\mu}_1'(I - P_1(\alpha))\boldsymbol{\mu}_1 = \boldsymbol{\mu}'B'(I - P_c)^{-1}B\boldsymbol{\mu},$$

where $P_c = X_{1c}(X'(\alpha)X(\alpha))^{-1}X_{1c}'$ (see, e.g., Result (5.4) of Chatterjee and Hadi [6], p. 189). One can now check that $(I - P_c)^{-1} = I + X_{1c}(X_1'X_1)^{-1}X_{1c}'$ and

$$(6.9) \qquad \boldsymbol{\mu}'B'(I - P_c)^{-1}B\boldsymbol{\mu} - \boldsymbol{\mu}'B'B\boldsymbol{\mu} = \boldsymbol{\mu}'B'X_{1c}(X_1'X_1)^{-1}X_{1c}'B\boldsymbol{\mu} \geq 0$$

as $(X_1'X_1)^{-1}$ is positive definite. From (6.7)–(6.9)

$$\frac{\boldsymbol{\mu}_1'(I - P_1(\alpha))\boldsymbol{\mu}_1}{nL_n(\alpha)} \leq \frac{\boldsymbol{\mu}_1'(I - P_1(\alpha))\boldsymbol{\mu}_1}{\boldsymbol{\mu}'(I - P(\alpha))\boldsymbol{\mu}} \leq \frac{||A\boldsymbol{\mu}||^2}{||A\boldsymbol{\mu}||^2 + ||B\boldsymbol{\mu}||^2}.$$

We now consider average over the $r$ training samples. Since each $y_i$ $(1 \leq i \leq n)$ appears in the same number of training samples, we have

$$(6.10) \qquad \frac{T_1}{nL_n(\alpha)} \leq \frac{\frac{1}{r}\sum_{i=1}^r \boldsymbol{\mu}_i'(I - P_i(\alpha))\boldsymbol{\mu}_i}{\boldsymbol{\mu}'(I - P(\alpha))\boldsymbol{\mu}} \leq \frac{n-k}{n}$$

which converges to zero.

To prove (6.6) for $j = 2$ we note that $T_2$ can be expressed as $\boldsymbol{e}'M(\alpha)\boldsymbol{e}$ for some matrix $M(\alpha) = (m_{ij})$, which is a sum of $r$ matrices corresponding to the $r$ choices of the $n - k$ indices from $\{1, 2, \ldots, n\}$ ($n - k$ rows of $X(\alpha)$). For example, for the training sample $\boldsymbol{y}_1 = (y_1, \ldots, y_{n-k})'$, $\boldsymbol{e}_1' P_1(\alpha) \boldsymbol{e}_1$ may be written as $\boldsymbol{e}' M_1(\alpha) \boldsymbol{e}$ where

$$M_1(\alpha) = \begin{pmatrix} P_1(\alpha) & \boldsymbol{0} \\ \boldsymbol{0} & \boldsymbol{0} \end{pmatrix}, \quad \text{thus} \quad M(\alpha) = (1/r) \sum_{i=1}^{r} M_i(\alpha).$$

As $P_i(\alpha)$'s are all idempotent matrices, one can show that $\sum_i \sum_j m_{ij}^2 \leq p_n(\alpha)$. Then proceeding as in the proof of (6.3) given in Li ([12], p. 970) one can prove the result using (6.2), (2.2), (2.3) and (6.5). Indeed, by (6.2),

$$P \left[ \max_\alpha \frac{|\boldsymbol{e}' M(\alpha) \boldsymbol{e} - \sigma^2 p_n(\alpha)|}{n R_n(\alpha)} > \epsilon \right]$$
$$\leq \quad C \sum_\alpha \frac{[p_n(\alpha)]^m}{\epsilon^{2m} [n R_n(\alpha)]^{2m}}$$

for some constant $C > 0$. The result follows from (2.2), (2.3) and (6.5).

The proof of (6.6) for $j = 3$ is similar. We note that $T_3 = \boldsymbol{e}' \boldsymbol{b}$ with $\boldsymbol{b} = (1/r) \sum_{i=1}^{r} (I - P_i(\alpha)) \boldsymbol{\mu}_i$ and $||b||^2 \leq (1/r) \sum_{i=1}^{r} \boldsymbol{\mu}_i' (I - P_i(\alpha)) \boldsymbol{\mu}_i$. By (6.1) and (6.10)

$$P \left[ \max_\alpha \frac{|\boldsymbol{e}' \boldsymbol{b}|}{n R_n(\alpha)} > \epsilon \right]$$
$$\leq \quad C \left( \frac{n-k}{n} \right)^m \sum_\alpha \frac{1}{[n R_n(\alpha)]^m}$$

for some constant $C > 0$. The result follows from (2.2) and (6.5). Thus (6.6) is proved and hence the lemma. $\qquad \square$

**Remark 6.1.** Indeed, to prove Lemma 2.1 and Lemma 2.2, we need to assume

$$\frac{p_n}{\min\limits_{\alpha \in \mathcal{A}_n} n R_n(\alpha)} \to 0$$

instead of the stronger condition (2.3). We, however, need (2.3) to prove our final result.

*Proof of Theorem 2.1.* Since (2.3) and (2.7) imply that the third term in the right hand side of (2.1) is of the order $o_p(L_n(\alpha))$ uniformly in $\alpha \in \mathcal{A}_n$, part (a) follows from (2.1), Lemma 2.1 and Lemma 2.2. From part (a), $\Gamma(\alpha)$ can be written as

$$\Gamma(\alpha) = \frac{k}{n^2} \boldsymbol{e}' \boldsymbol{e} + L_n(\alpha)(1 + \zeta_n(\alpha)), \ \ \alpha \in \mathcal{A}_n,$$

where $\max\limits_\alpha |\zeta_n(\alpha)| \xrightarrow{p} 0$. Now $\Gamma(\hat{\alpha}_n) \leq \Gamma(\alpha) \ \forall \ \alpha$ implies

$$\frac{L_n(\hat{\alpha}_n)}{L_n(\alpha)} \leq \frac{1 + \zeta_n(\alpha)}{1 + \zeta_n(\hat{\alpha}_n)} \leq \frac{1 + \max\limits_\alpha |\zeta_n(\alpha)|}{1 - \max\limits_\alpha |\zeta_n(\alpha)|} \quad \forall \alpha.$$

Part (b) follows from the above. $\qquad \square$

*Proof of Lemma 3.1.* We first note that under suitable conditions there exist $0 < \delta < \Delta$ such that

$$(6.11) \qquad \qquad \log(1 + \delta) < u_n(\alpha) < \log(1 + \Delta) \quad \forall \alpha$$

with probability tending to 1. This follows from (3.3), (3.4) and the fact that $\boldsymbol{e}'\boldsymbol{e}/n\sigma^2 \xrightarrow{p} 1$, noting that $\max_\alpha \boldsymbol{e}'P(\alpha)\boldsymbol{e}/n \leq \boldsymbol{e}'P\boldsymbol{e}/n \xrightarrow{p} 0$ and $L_n(\alpha)$ is uniformly (in $\alpha$) bounded with probability tending to 1. Here $P$ is the projection matrix corresponding to the full model.

Consider now the expression in (3.2). By Lemma 2.1 of Section 2 and (6.11),

$$
\begin{aligned}
\log[S(\alpha)/n\sigma^2] &= \log[\boldsymbol{e}'\boldsymbol{e}/n\sigma^2 + L_n(\alpha)/\sigma^2 + o_p(L_n(\alpha)/\sigma^2)] \\
&= \log[\boldsymbol{e}'\boldsymbol{e}/n\sigma^2 + L_n(\alpha)/\sigma^2 + o_p(\boldsymbol{e}'\boldsymbol{e}/n\sigma^2 + L_n(\alpha)/\sigma^2)] \\
&= \log[e^{u_n(\alpha)}(1 + o_p(1))] \\
&= u_n(\alpha) + o_p(1) \\
(6.12) \qquad &= u_n(\alpha) + o_p(u_n(\alpha))
\end{aligned}
$$

uniformly in $\alpha$. In view of (3.5), to prove (3.7), it remains to show

$$(6.13) \qquad \qquad \frac{n-k}{nr} \sum_{i=1}^r \log[S_i(\alpha)/n\sigma^2] = o_p(1).$$

Note that we are also using (3.4) and (6.11). Since $S_i(\alpha) \geq S_i$ for all $\alpha$ and all $i$, we have for all $\alpha$

$$0 < \frac{1}{r} \sum_{i=1}^r \log[S_i(\alpha)] = \log[\prod_{i=1}^r S_i(\alpha)]^{1/r} \leq \log[\frac{1}{r} \sum_{i=1}^r S_i(\alpha)]$$

implying

$$
\begin{aligned}
-\frac{n-k}{n}\log(n\sigma^2) &< \frac{n-k}{nr} \sum_{i=1}^r \log\left[\frac{S_i(\alpha)}{n\sigma^2}\right] \\
&\leq \frac{n-k}{n}\log[\frac{1}{r}\sum_{i=1}^r S_i(\alpha)] - \frac{n-k}{n}\log(n\sigma^2).
\end{aligned}
$$

Then (6.13) follows from Lemma 2.2 of Section 2, condition (3.4) and the fact that $L_n(\alpha)$ is uniformly (in $\alpha$) bounded with probability tending to 1 (as noted earlier in the argument for (6.11)). $\qquad \square$

*Proof of Theorem 3.1.* Let $\hat{\alpha}_n$ be the model which minimizes $\Gamma(\alpha)$. Proceeding as in the proof of part (b) of Theorem 2.1, and using (3.7) we can prove that

$$\frac{u_n(\hat{\alpha}_n)}{u_n(\alpha_n^L)} \xrightarrow{p} 1.$$

This, together with (6.11), imply that

$$u_n(\hat{\alpha}_n) - u_n(\alpha_n^L) \xrightarrow{p} 0$$

i.e., $\qquad \dfrac{\boldsymbol{e}'\boldsymbol{e} + nL_n(\hat{\alpha}_n)}{\boldsymbol{e}'\boldsymbol{e} + nL_n(\alpha_n^L)} \xrightarrow{p} 1.$

Since $\frac{\boldsymbol{e'e}}{n} \xrightarrow{p} \sigma^2$ and $L_n(\alpha_n^L) \geq \min_\alpha \Delta_n(\alpha)$, using (3.3) we have

$$\frac{L_n(\hat{\alpha}_n)}{L_n(\alpha_n^L)} \xrightarrow{p} 1.$$

$\square$

*Proof of Proposition 4.1.* We first prove equation (4.7). Below, by $\max\limits_\alpha$ we mean maximum over $\alpha \in \mathcal{A}_n^c$. Let $Z_n(\alpha) = (\boldsymbol{e'}P(\alpha)\boldsymbol{e})/(\sigma^2 p_n(\alpha))$. We first show that $\max\limits_\alpha |Z_n(\alpha)| = O_p(1)$. By (6.2)

$$P[\max_\alpha |Z_n(\alpha) - 1| > M]$$
$$\leq \sum_\alpha E|Z_n(\alpha) - 1|^{2m}/M^{2m}$$
$$\leq \frac{C}{M^{2m}} \sum_\alpha \frac{1}{[p_n(\alpha)]^m}$$

for some constant $C > 0$ and by (4.5) this can be made arbitrarily small by choosing suitable $M > 0$. Thus $\max\limits_\alpha |Z_n(\alpha) - 1| = O_p(1)$ implying $\max\limits_\alpha |Z_n(\alpha)| = O_p(1)$. This implies $(1/n)\boldsymbol{e'}P(\alpha)\boldsymbol{e} = o_p(\frac{1}{n}\lambda_n\sigma^2 p_n(\alpha))$ uniformly in $\alpha \in \mathcal{A}_n^c$ as $\lambda_n \to \infty$. Proceeding in a similar manner and noting that $(1/r)\sum\limits_{i=1}^r \boldsymbol{e}_i'P_i(\alpha)\boldsymbol{e}_i$ can be written as $\boldsymbol{e'}M(\alpha)\boldsymbol{e}$ (see proof of Lemma 2.2) one can prove

$$\frac{1}{nr}\sum_{i=1}^r \boldsymbol{e}_i'P_i(\alpha)\boldsymbol{e}_i = o_p(\frac{1}{n}\lambda_n\sigma^2 p_n(\alpha)) \text{ uniformly in } \alpha \in \mathcal{A}_n^c.$$

The result now follows from (4.2), (4.4) and (4.6).
In order to complete the proof of Proposition 4.1, we now prove equation (4.8). From (4.7),

$$\Gamma(\alpha_n^c) = \frac{k}{n^2}\boldsymbol{e'e} + \frac{1}{n}\lambda_n\sigma^2 p_n(\alpha_n^c) + o_p\left(\frac{1}{n}\lambda_n\sigma^2 p_n(\alpha_n^c)\right).$$

The result follows from (4.1) and (4.2) noting that (4.1) implies (6.5) with $\max\limits_{\alpha \in \mathcal{A}_n}$ replaced by $\max\limits_{\alpha \in \mathcal{A}_n - \mathcal{A}_n^c}$. $\square$

*Proof of (4.18).* Note that

$$\max_{\alpha \in \mathcal{A}_n - \mathcal{A}_n^c} \frac{L_n(\alpha_n^c)}{L_n(\alpha)} = \max_\alpha \frac{\boldsymbol{e'}P(\alpha_n^c)\boldsymbol{e}}{nL_n(\alpha)}.$$

By (6.2) and by arguments used earlier

$$P\left[\max_{\alpha \in \mathcal{A}_n - \mathcal{A}_n^c} \left|\frac{\boldsymbol{e'}P(\alpha_n^c)\boldsymbol{e} - \sigma^2 p_n(\alpha_n^c)}{nR_n(\alpha)}\right| > \epsilon\right]$$
$$\leq C\left[\frac{p_n}{\min\limits_\alpha nR_n(\alpha)}\right]^m \sum_{\alpha \in \mathcal{A}_n - \mathcal{A}_n^c} \frac{1}{[nR_n(\alpha)]^m}$$

for some constant $C$. The result follows from (4.1) and (4.2). $\square$

**Acknowledgments.** We would first like to thank Professors Bertrand Clarke and Subhashis Ghosal for inviting us to contribute to this volume. We consider it a great privilege to be able to write an article as an expression of our deep regard for our mentor and teacher Professor Jayanta K. Ghosh. We cannot express in words how much we learned in numerous academic discussions we had with him over the years which immensely influenced our thought process. Needless to say, this work owes much to the way we learned to think about model selection through our association with him. Finally, we thank an anonymous referee and the editors for helpful comments and suggestions towards improvement of the paper.